\documentclass[graybox]{svmult}
\usepackage{mathptmx}       
\usepackage{helvet}         
\usepackage{courier}        
\usepackage{type1cm}        
\usepackage{makeidx}         
\usepackage{graphicx}        
\usepackage{multicol}        
\usepackage[bottom]{footmisc}
\usepackage{tikz} 

\usepackage{url}

\usepackage{amsfonts}
\usepackage{amsmath}
\usepackage{amssymb}

\newcommand{\R}{\mathbb{R}}

\newcommand{\pd}[2]{\frac{\partial #1}{\partial #2}}

\newcommand{\td}[2]{\frac{\mathrm d #1}{\mathrm d #2}}
  
\newcommand{\dx}{\,\mathrm{d}\mathbf{x}}
\newcommand{\ds}{\,\mathrm{d}s}

\newcommand{\beq}{\begin{equation}}
\newcommand{\eeq}{\end{equation}}

\newcommand{\Ne}{\mathcal N^e}
\newcommand{\Ei}{\mathcal E_i}

\begin{document}

\title*{Convex limiting for finite elements and its relationship to residual distribution}
\author{Dmitri Kuzmin} 
\institute{Prof. Dmitri Kuzmin \at
  Institute of Applied Mathematics (LS III), TU Dortmund University\\
  Vogelpothsweg 87, D-44227 Dortmund, Germany
  \at \email{kuzmin@math.tu-dortmund.de}}

\maketitle

\abstract{We review some recent advances in the field of element-based
algebraic stabilization for continuous finite element discretizations
of nonlinear hyperbolic problems. The main focus is on multidimensional
convex limiting techniques designed to constrain antidiffusive element
contributions rather than fluxes. We show that the resulting 
schemes can be interpreted as residual distribution methods.
Two kinds of convex limiting can be used to enforce the validity of
generalized discrete maximum principles in this context. The first
approach has the structure of a localized flux-corrected transport (FCT)
algorithm, in which the computation of~a low-order predictor is followed by
an antidiffusive correction stage. The second option is the use of a
monolithic convex limiting (MCL) procedure at the level of
spatial semi-discretization. In both cases, inequality constraints
are imposed on scalar functions of intermediate states that are
required to stay in convex invariant sets. 
}




\section{Introduction}
\label{sec:1}

The design of a finite element scheme that ensures preservation of local
and/or global bounds for numerical solutions to a hyperbolic problem
typically involves construction and local adjustment of a dissipative
shock capturing term. In most cases, a finite element counterpart of
the local Lax--Friedrichs (Rusanov) method is corrected by adding
limited antidiffusive terms. In edge-based schemes, the additional
terms consist of fluxes between pairs of neighboring nodes
\cite{guermond2018,convex,fctbook,kuzmin2010,luo1994}. In
element-based alternatives, diffusive and antidiffusive
correction terms are assembled from element vectors that
possess the zero sum property \cite{dobrev2018,preprint,kuzmin2020,lohmann2017}.
A~comprehensive review of the state of the art can be found in 
the recent book \cite{book}.

The element-based approach traces
its origins to the flux-corrected transport method
proposed by Löhner et al. \cite{lohner1987,lohner1988}. Moreover,
it is closely related to residual distribution (RD) schemes that
achieve the algebraic stabilization effect by manipulating
element contributions in a locally conservative manner
\cite{abgrall2006,abgrall2017,abgrall2019,abgrall2010}. In this
note, we review the recent progress in the development of
element-based convex limiting / RD tools for continuous Galerkin
methods using linear finite elements.

\section{Element-based Rusanov method}
\label{sec:2}

Let $u(\mathbf x,t)$ denote an exact solution of the
hyperbolic conservation law or system
\beq\label{conslaw}
 \pd{u}{t}+\nabla\cdot\mathbf{f}(u)=0\quad
 \mbox{in}\ \Omega\times (0,T].
   \eeq
We assume that periodic boundary conditions are imposed
on $\partial\Omega$ and that there exists a convex invariant set
$\mathcal G$ such that $u\in \mathcal G$ a.e. in $\Omega\times (0,T]$. 

A continuous piecewise-linear ($\mathbb P_1$)
finite element approximation $u_h=\sum_{j=1}^{N_h}u_j\varphi_j$ to $u$ is called
\emph{invariant domain preserving} (IDP) or positivity preserving if
$u_j\in \mathcal G$ for $j=1,\dots,N_h$. The Lagrange basis functions
$\varphi_i$ are associated with vertices $\mathbf x_i$ of a conforming mesh
$\mathcal T_h$. The elements (cells) of that mesh are denoted by
$K_1,\ldots,K_{E_h}$. The global indices of nodal points $\mathbf x_i$
belonging to $K_e$ are stored in the set $\Ne$. The
indices of elements $K_e$ containing  $\mathbf x_i$ are stored
in the set $\Ei$.

The standard Galerkin discretization of the hyperbolic problem
\eqref{conslaw} is given by
\beq\label{galerkin}
\sum_{e\in\mathcal E_i}m_{ij}^e\td{u_j}{t}=\sum_{e\in\mathcal E_i}
\int_{K_e}\nabla\varphi_i\cdot\mathbf f(u_h)\dx,\qquad i=1,\ldots,N_h,
\eeq
where $m_{ij}^e=\int_{K_e}\varphi_i\varphi_j\dx$ is the contribution of
$K_e$ to the entry $m_{ij}=\sum_{e\in\mathcal E_i\cap\mathcal E_j}m_{ij}^e$ of the
consistent mass matrix $M_C=\{m_{ij}\}_{i,j=1}^{N_h}$. 

To construct a low-order IDP approximation to \eqref{galerkin}, we lump
the mass matrix, use the midpoint rule for the volume integrals on
the right-hand side, and add Rusanov artificial
viscosity of the form employed in \cite{abgrall2006}. These
manipulations lead to
\beq\label{rusanov}
m_i\td{u_i}{t}=\sum_{e\in\mathcal E_i}
[d^e(\bar u^e-u_i)-\mathbf f(\bar u^e)
  \cdot\mathbf c_i^e],\qquad i=1,\ldots,N_h,
\eeq
where $\bar u_e$ is the arithmetic mean of $u_i,\ i\in\Ne$ and
$\mathbf  c_i^e=-\int_{K_e}\nabla\varphi_i\dx$. The Rusanov
viscosity coefficient $d^e=\max_{i\in\Ne}\lambda_i^e|\mathbf  c_i^e|$
is defined using the maximum speeds $\lambda_i^e$ of 
Riemann problems with flux functions
$\mathbf f(u)\cdot \mathbf  c_i^e/|\mathbf  c_i^e|$
and initial states $(\bar u^e,u_i)$.

If \eqref{rusanov}
is discretized in time using an explicit strong stability
preserving (SSP) Runge--Kutta method, then each low-order forward Euler stage
can be written as 
\beq\label{rusanovSSP}
u_i^L=\Big(1-\frac{2\Delta t}{m_i}\sum_{e\in\mathcal E_i}d^e\Big)u_i
+\frac{2\Delta t}{m_i}\sum_{e\in\mathcal E_i}d^e\bar u_i^e.
\eeq
Similarly to the \emph{bar states} of the edge-based version analyzed
in \cite{guermond2016}, the state
$$
\bar u_i^e=\frac{\bar u^e+u_i}2-
\frac{(\mathbf f(\bar u^e)-\mathbf f(u_i))\cdot\mathbf c_i^e}{2d_i^e},\qquad
e\in\Ei
$$
represents a spatially
averaged exact solution of the associated Riemann problem. Thus,
$\bar u^e,u_i\in\mathcal G$ implies $\bar u_i^e\in\mathcal G$ for
any convex invariant set $\mathcal G$. If the time step $\Delta t$
satisfies the CFL-like condition
$\frac{2\Delta t}{m_i}\sum_{e\in\mathcal E_i}d^e\le 1$, then $u_i^L$
is a convex combination of the IDP states $u_i$ and $\bar u_i^e,\
e\in\mathcal E_i$. Hence, the fully discrete scheme is IDP.

The element-based low-order method \eqref{rusanov} is new but
represents a slightly modified version of the one derived
in \cite[Sec. 4.6.1]{book} using
a trapezoidal rule approximation for the flux
divergence term. The use of the midpoint rule
reduces the effort for calculating the Rusanov viscosity
$d^e$ and makes it possible
to prove the IDP property without introducing edge-based bar
states (cf. \cite[Thm. 4.13]{book}).

\section{Residual distribution}
\label{sec:3}

To interpret \eqref{rusanov} as a low-order RD method derived from 
\eqref{galerkin}, we notice that both semi-discrete schemes can
be written in the generic form (cf. \cite[Chapter 4]{book})
\beq\label{rdform}
m_i\td{u_i}{t}=\sum_{e\in\mathcal E_i}r_i^e,
\eeq
where $r_i^e$ denotes the contribution of $K_e$ to the steady-state
residual. Using integration by parts, we deduce that \eqref{galerkin}
is equivalent to \eqref{rdform} for $r_i^e=r_i^{e,H}$ given by
$$
r_i^{e,H}=\int_{K_e}\varphi_i(\dot u_i-\dot u_h)\dx
-\int_{K_e}\varphi_i\nabla\cdot\mathbf f(u_h)\dx.
$$
The coefficients $\dot u_i=\td{u_i}{t}$ of 
$\dot u_h=\td{u_h}{t}$ can be calculated using a matrix-free
iterative solver for system
\eqref{galerkin} or a truncated Neumann series
approximation to $M_C^{-1}$. Such algorithms for inversion of
mass matrices can be found in
\cite{abgrall2017,donea1984,guermond2014}. For linear finite
elements, which we are using in this work, it is worthwhile to approximate $\dot u_i$ by
$$
\dot u_i^L=\frac{1}{m_i}
\sum_{e\in\mathcal E_i}
    [d^e(\bar u^e-u_i)-\mathbf f(\bar u^e)\cdot\mathbf c_i^e]
=\frac{1}{m_i}
\sum_{e\in\mathcal E_i}2d^e(\bar u_i^e-u_i).
$$
This approximation corresponds to adding a
high-order stabilization term \cite{convex}.

In the low-order method \eqref{rusanov}, the element contribution
$r_i^{e,H}$ is replaced with
$$
r_i^{e,L}=r_i^{e,H}-f_i^e,
$$
$$
f_i^e=\int_{K_e}\varphi_i(\dot u_i-\dot u_h)\dx
+\int_{K_e}\nabla\varphi_i\cdot(\mathbf f(u_h)-\mathbf f(\bar u^e))\dx
-d^e(\bar u^e-u_i).
$$
The Lagrange basis functions $\varphi_i$ have the property that
$\sum_{i\in\Ne}\varphi_i\equiv 1$. Hence,
$\sum_{i\in\Ne}\nabla \varphi_i\equiv \mathbf 0$ and
$\sum_{i\in\Ne}
\int_{K_e}\varphi_i(\dot u_i-\dot u_h)\dx=0$.
Furthermore, $\sum_{i\in\Ne}(\bar u^e-u_i)=0$
by definition of $\bar u^e$. It follows that the element
contributions satisfy
$$
\sum_{i\in\Ne}f_i^e=0,\qquad \sum_{i\in\Ne}r_i^{e,L}
=\sum_{i\in\Ne}r_i^{e,H}=r^e,
$$
\beq\label{fluctu}
r^e=-\int_{K_e}\nabla\cdot\mathbf f(u_h)\dx
=-\int_{\partial K_e}\mathbf f(u_h)\cdot\mathbf n\ds,
\eeq
where $\mathbf n$ is the unit outward normal.
The general form \eqref{rdform} with
$
r_i^{e}=r_i^{e,L}+f_i^{e,*}
$
reduces to \eqref{galerkin} for $f_i^{e,*}=f_i^e$ and to
\eqref{rusanov} for $f_i^{e,*}=0$. The \emph{fluctuation}
$r^e=\sum_{i\in\Ne}r_i^{e}$ defined in formula \eqref{fluctu}
is preserved provided that
$\sum_{i\in\Ne}f_i^{e,*}=0$.

A classical RD method for a scalar conservation law of the form 
\eqref{conslaw} constructs nodal residuals $r_i^e=\beta_i^er^e$
using distribution weights $\beta_i^e$ such that $\sum_{i\in\Ne}
\beta_i^e=1$. If $r^e\ne 0$ and $r_i^e$ is given, then
$\beta_i^e=r_i^e/r^e$. For $r^e=0$, the assumption
that $r_i^e=\beta_i^er^e$ implies $r_i^e=0\ \forall i\in\Ne$.
However, $r^e=0$ holds also if $r^e_+=-r^e_-\ne 0$, where
$$
r^e_+=\sum_{i\in\Ne}\max\{0,r_i^e\},\qquad
r^e_-=\sum_{i\in\Ne}\min\{0,r_i^e\}.
$$
The authors of \cite{RD1,RD2} distribute
the fluctuations
$r^e_\pm$ using weights $\beta_{i,\pm}^{e}$ such that 
$$
r_i^e=\begin{cases}
\beta_{i,+}^er^e_+ & \mbox{if}\ r_i^e>0,\\
\beta_{i,-}^er^e_- & \mbox{if}\ r_i^e<0,
\end{cases}\qquad \beta_{i,\pm}^{e}\ge 0,\quad
\sum_{i\in\Ne}\beta_{i,\pm}^e=1.
$$
As we show in the next section, such distribution can also be used
for $(f_i^{e,*})_{i\in\Ne}$.

\section{Convex limiting}
\label{sec:4}

Recall that the IDP property of the low-order scheme
was shown by representing the result of the explicit update
\eqref{rusanovSSP} as a convex combination of states
belonging to a convex invariant set $\mathcal G$. Modern
convex limiting techniques ensure the existence of such
representations for high-order extensions. Element-based
limiters of this kind can be subdivided into (A)
flux-corrected transport (FCT) algorithms that add
$f_i^{e,*}$ to a low-order predictor \cite{RD1,book,lohmann2017}
and (B) monolithic semi-discrete alternatives that insert $f_i^{e,*}$
into the right-hand side of the low-order scheme \cite{RD2,book,preprint}.

In the first stage of a type A scheme, a low-order IDP
approximation $u_h^L$ is calculated using \eqref{rusanovSSP}.
The second stage corrects $u_i^L\in\mathcal G$ as follows:
$$
u_i^{\rm CL}=u_i^L+\frac{\Delta t}{m_i}\sum_{e\in\Ei}f_i^{e,*}=
\frac{1}{m_i}\sum_{e\in\Ei}m_i^e\bar u_i^{e,*},\qquad
\bar u_i^{e,*}=u_i^L+\frac{\Delta tf_i^{e,*}}{m_i^e}.
$$
Since $m_i=\sum_{e\in\Ei}m_i^e$ with $m_i^e>0$,
the result $u_i^{\rm CL}$ is a convex
combination of the auxiliary states $\bar u_i^{e,*}$. The element-based
limiter should guarantee that $\bar u_i^{e,*}\in\mathcal G$.
The first representative of such limiters was introduced in
\cite[Sec. 2.3.2]{cotter2016} in the context of flux-corrected
remapping for scalar quantities. Element-based FCT/RD algorithms
using this limiting strategy can be found in
\cite{anderson2017,dobrev2018,RD1,kuzmin2020,lohmann2017}.
An edge-based extension to nonlinear hyperbolic systems was proposed
by Guermond et al.~\cite{guermond2018}.

A monolithic convex limiting (MCL) algorithm
is a single-stage type B scheme that avoids the computation
of $u_i^L$ and replaces \eqref{rusanovSSP} with (cf. \cite{convex})
$$
u_i^{\rm CL}=\Big(1-\frac{2\Delta t}{m_i}\sum_{e\in\mathcal E_i}d^e\Big)u_i
+\frac{2\Delta t}{m_i}\sum_{e\in\mathcal E_i}d^e\bar u_i^{e,*},\qquad
 \bar u_i^{e,*}=\bar u_i^e+\frac{f_i^{e,*}}{2d^e}.
$$
The IDP property is guaranteed under the CFL condition
of the low-order scheme if $\bar u_i^{e,*}\in\mathcal G\
\forall e\in\mathcal E_i$. Representatives of element-based
MCL schemes can be found in \cite[Chap. 4, 6]{book}.
The original edge-based version was introduced in \cite{convex}.

Limiters of type A and B have a lot in common.
In both versions, $\bar u_i^{e,*}$ consists of a low-order IDP
state $\bar u_i^{e,L}\in\{u_i^L,\bar u_i^e\}$ and a correction term
$f_i^{e,*}/\gamma_i^e$ with $\gamma_i^e\in\{m_i^e/\Delta t,2d^e\}$.
Hence, the same algorithms can be used to enforce inequality
constraints for $\bar u_i^{e,*}=\bar u_i^{e,L}+f_i^{e,*}/\gamma_i^e$
in type A and type B schemes \cite{convex,book}.

\subsection{Limiting for scalars}

In the scalar case, the states $\bar u_i^{e,*}$
can be constrained to stay in a local admissible range
$[u_i^{\min},u_i^{\max}]\subseteq[u^{\min},u^{\max}]=:\mathcal G$
such that $\bar u_i^{e,L}\in[u_i^{\min},u_i^{\max}]\ \forall e\in\Ei$.
The constraints for element-based convex limiting of type A or B
are formulated as follows:
\begin{gather}
f_{i}^{e,\min}:=
\gamma_i^e(u_i^{\min}-\bar u_i^{e,L})\le f_{i}^{e,*}
\le
\gamma_i^e(u_i^{\max}-\bar u_i^{e,L})=: f_{i}^{e,\max},\label{dmp}\\
\sum_{i\in\Ne}f_{i}^{e,*}=0.\label{zerosum}
\end{gather}
The zero sum condition \eqref{zerosum} is clearly satisfied for
$f_{i}^{e,*}=\alpha^ef_{i}^e$, where $\alpha^e\in[0,1]$ is
an element-based correction factor. Conditions \eqref{dmp}
are satisfied for \cite{cotter2016}
\beq\label{alphae}
\alpha^e=\min_{i\in\Ne}\alpha_{i}^e,\qquad
\alpha_{i}^e=
\begin{cases}
f_{i}^{e,\max}/f_{i}^e &
\mbox{if}\ \ f_{i}^e>f_{i}^{e,\max},\\
f_{i}^{e,\min}/f_{i}^e &
\mbox{if}\ \ f_{i}^e<f_{i}^{e,\min},\\
1 & \mbox{otherwise}.
\end{cases}
\eeq
This \emph{scaling} limiter was used in \cite{anderson2017,RD2,lohmann2017}
to constrain finite element approximations of very high order. Interestingly
enough, \eqref{alphae} has the same structure as the Barth--Jespersen slope
limiter \cite{barth1989} for unstructured grid finite volume methods.

The application of $\alpha^e$ to $f_i^{e}=\beta_{i,\pm}^ef^e_\pm$ can be
interpreted as residual distribution that scales the fluctuations
$f^e_\pm$, while leaving the weights $\beta_{i,\pm}^e$ unchanged.

Another way to satisfy the limiting constraints \eqref{dmp} and \eqref{zerosum}
for each node $i\in\Ne$ is the following \emph{clip-and-scale}
(C\&S) limiting strategy
\cite{anderson2017,RD1,book,lohmann2017}:
\begin{enumerate}
\item Calculate
$\tilde f_i^e=\max\{f_i^{e,\min},
\min\{f_i^e,f_i^{e,\max}\}\}
$
and the fluctuations
\[
\tilde f^e_{+}=\sum_{i\in\Ne}\max\{0,\tilde f_i^e\},\qquad
\tilde f^e_{-}=\sum_{i\in\Ne}\min\{0,\tilde f_i^e\}.
\]
\item Balance the sums of positive and negative components 
\beq\label{csfel}
f_i^{e,*}=\begin{cases}
  - (\tilde f^e_-/\tilde f^e_+)\tilde f_i^e &
\mbox{if}\ \tilde f_i^e>0,\ 
\tilde f^e_++\tilde f^e_->0,\\
- (\tilde f^e_+/\tilde f^e_-)\tilde f_i^e &
\mbox{if}\ \tilde f_i^e<0, \ 
\tilde f^e_++\tilde f^e_-<0,\\
\tilde f_i^e & \mbox{otherwise}.\end{cases}
\eeq
\end{enumerate}
The RD interpretation of this element-based limiter is as follows. The first
step replaces $f_i^{e}=\beta_{i,\pm}^ef^e_\pm$ with clipped antidiffusive
element contributions $\tilde f_i^e=\tilde\beta_{i,\pm}^e\tilde f^e_\pm$
that satisfy \eqref{dmp}. The second step enforces the zero sum
condition \eqref{zerosum} by scaling the positive or negative
fluctuation. The weights $\tilde\beta_{i,\pm}^e$
remain unchanged.

Note that $\tilde f_i^e=\tilde\alpha_i^ef_i^e$, where
$\tilde\alpha_i^e$ corresponds to $\alpha_i^e$
defined by \eqref{alphae}.
Hence, the result of C\&S limiting 
can be written as $f_i^{e,*}=\alpha_i^ef_i^e$, where
$\alpha_i^e=-(\tilde f^e_\mp/\tilde f^e_\pm) \tilde\alpha_i^e$ if
$\tilde f^e_++\tilde f^e_-\ne 0$.
That is, the C\&S algorithm applies an individually chosen correction
factor $\alpha_i^e$ to each component $f_i^e$. As shown in \cite{lohmann2017},
this nodal limiting
strategy is less diffusive than simple scaling using the common
correction factor \eqref{alphae}. Moreover, the limited element
contributions \eqref{csfel} depend continuously on the data.

Further examples of scalar element-based limiting techniques
can be found in \cite[Chap. 4]{book}.
A classical representative of RD-based nodal limiters
is the nonlinear positive streamwise invariant
(PSI) method \cite{abgrall2006,RD1,struijs1994}. Such limiters
produce distribution
weights $\beta_{i,\pm}^e$ instead of correction factors $\alpha_i^e$.
We remark that, while representations of $f_i^{e,*}$ in terms of 
 $\alpha_i^e$ or $\beta_{i,\pm}^e$ may be useful for
derivation and comparison purposes, 
direct computation of $f_i^{e,*}$ should be preferred in practice.

\subsection{Limiting for systems}

Let us now consider a system of conservation
laws for
$u=(\varrho,\varrho\phi_2,\ldots,\varrho\phi_{m})^\top$,
where $\varrho=\rho\phi_1$ with $\phi_1=1$ is a density-like variable.
The invariant set to be preserved by an IDP limiter for the element
contributions
$f_i^e=(f_{i,\varrho}^{e},f_{i,\varrho\phi_1}^{e},\ldots,f_{i,\varrho\phi_m}^{e})^\top$
is usually of the form
$\mathcal G=\{u\in\R^m\,:\,\Phi_1(u)\ge0,\ldots,\Phi_{M}(u)\ge 0\}$,
where  $\Phi_1,\ldots,\Phi_M$ are quasi-concave scalar functions
of the conserved variables. The quasi-concavity implies that
$\Phi(u)\ge 0$ for any convex combination of intermediate states
$\bar u$ such that $\Phi(\bar u)\ge 0$. This property is exploited in
low-order IDP schemes and in limiting techniques based on convex
decompositions \cite{guermond2016,convex,book,wu2025}.
  
In sequential convex limiting algorithms \cite{dobrev2018,book,kuzmin2020},
local maximum principles for $\rho,\phi_2,\ldots,\phi_m$ 
and the IDP property w.r.t. $\mathcal G$ are enforced as follows:
\begin{enumerate}
\item Given  $f_{i,\varrho}^{e}$, calculate $ f_{i,\varrho}^{e,*}$
using a scalar (scaling or C\&S) limiter to ensure that
\[
\varrho_i^{\min}
\le \bar\varrho_i^{e,*}=\bar\varrho_i^{e,L}+f_{i,\varrho}^{e,*}/\gamma_i^e\le\varrho_i^{\max}.\]

\item For $\phi\in\{\phi_2,\ldots,\phi_m\}$, use a product rule version
  of the scalar limiter (see below) to construct element contributions
  $ f_{i,\varrho\phi}^{e,*}\approx
  f_{i,\varrho\phi}^{e}$ satisfying
    \[
    \bar\varrho_i^{e,*}\phi_i^{\min}\le \bar\varrho_i^{e,*}\bar\phi_i^{e,*}=
    \overline{(\varrho\phi)}_i^{e,L}+f_{i,\varrho\phi}^{e,*}/\gamma_i^e
    \le\bar\varrho_i^{e,*}\phi_i^{\max}.
\]
\item To keep $\bar u_i^{e,*}$ in the convex admissible set $\mathcal G$, apply
  $\alpha^e\in[0,1]$ such that
  \[\bar u_i^e\in\mathcal G\quad\Rightarrow\quad
  \bar u_i^{e,*}=\bar u_i^{e,L}+\alpha^e f_i^{e,*}/\gamma_i^e
  \in\mathcal G.
  \]
\end{enumerate}

In Step 2, the local maximum principle 
$\phi_i^{\min}\le\bar\phi_i^{e,*}\le\phi_i^{\max}$ is enforced using a discrete
version of $(\rho\phi)'=\rho'\phi+\rho\phi'$.
Simple scaling limiters designed for this purpose can
be found in \cite{dobrev2018,kuzmin2020}. A product rule
version of the C\&S limiter was introduced in
\cite[Sec. 6.3.2]{book}. It adjusts $f_{\varrho\phi}^{e}
=(f_{i,\varrho\phi}^{e})_{i\in\Ne}$
in the following way:
\begin{center}
\begin{minipage}{0.9\textwidth}
  \begin{itemize}
    
  \item Construct
    $\delta f_{\varrho\phi}^{e}=(\bar\phi_i^{e} f_{i,\varrho}^{e,*})_{i\in\Ne}$
    using the nodal states
    $\bar\phi_i^{e}=
    \frac{\overline{(\varrho\phi)}_i^{e,L}}{\bar \varrho_i^{e,L}}$.
    
\item Apply the scaling operator $\mathcal R_S$ to
  $\delta f^e_{\varrho\phi}$ and calculate
    $$
    g^{e}_{\varrho\phi}=f^{e}_{\varrho\phi}-\mathcal R_S(\delta f^{e}_{\varrho\phi}).
    $$
    
  \item Construct $\phi_i^{\min}=\min_{e\in\mathcal E_i}\bar\phi_i^{e,L}$
     and $\phi_i^{\max}=\max_{e\in\mathcal E_i}\bar\phi_i^{e,L}$ using
    $$
    \bar\phi_i^{e,L}= \frac{1}{\bar\varrho_i^{e,*}}
     \Big[  \overline{(\varrho\phi)}_i^{e,L}+
       (f_{i,\varrho\phi}^{e}-g_{i,\varrho\phi}^{e})/\gamma_i^e
       \Big].
     $$
     
   \item Use the scalar C\&S limiter to calculate
     $g_{i,\varrho\phi}^{e,*}$ such that
     $$
       \bar\varrho_i^{e,*}\phi_i^{\min}\le
       \bar\rho_i^{e,*}\bar\phi_i^{e,L}+g_{i,\varrho\phi}^{e,*}/\gamma_i^e
       \le \bar\varrho_i^{e,*}\phi_i^{\max},\qquad
  \sum_{i\in\Ne}g_{i,\varrho\phi}^{e,*}=0.
   $$
   \item Calculate 
     $
     f_{i,\varrho\phi}^{e,*}=f_{i,\varrho\phi}^{e}-g_{i,\varrho\phi}^{e}
     +g_{i,\varrho\phi}^{e,*}.
     $
  \end{itemize}
\end{minipage}
\end{center}
If $\Phi_1=\varrho$ is the only variable that should stay
nonnegative by definition of $\mathcal G$, then Step 3
of the sequential limiting procedure will
use $\alpha^e=1$, because the local maximum principle
$0\le\varrho_i^{\min}\le\bar\varrho_i^{e,*}\le\varrho_i^{\max}$ is
enforced in Step 1. Each additional variable
$\Phi\in\{\Phi_2,\ldots,\Phi_M\}$ imposes an upper bound
$\alpha_\Phi^e$ on the value of $\alpha^e$. This bound
can be calculated using a general line search algorithm
\cite{guermond2018} or a closed-form expression derived
from linear sufficient conditions \cite{convex,book,wu2023}.
Limiters that ensure positivity preservation for the
pressure and internal energy of the compressible Euler
equations in this way can be found in \cite{convex,kuzmin2020,abgrall2025}
and \cite[Example 4.16]{book}.

In principle, it is possible to preserve local and/or global bounds
for all scalar quantities of interest using the same scaling
factor $\alpha^e$ for all components of $f_i^e$. However, such
synchronized limiting is likely to produce more diffusive results.

\subsection{Example: Euler equations}

\begin{figure}[h!]

  \centering

  (a) low-order solution

\includegraphics[width=0.8\textwidth,trim=50 300 250 290,clip]{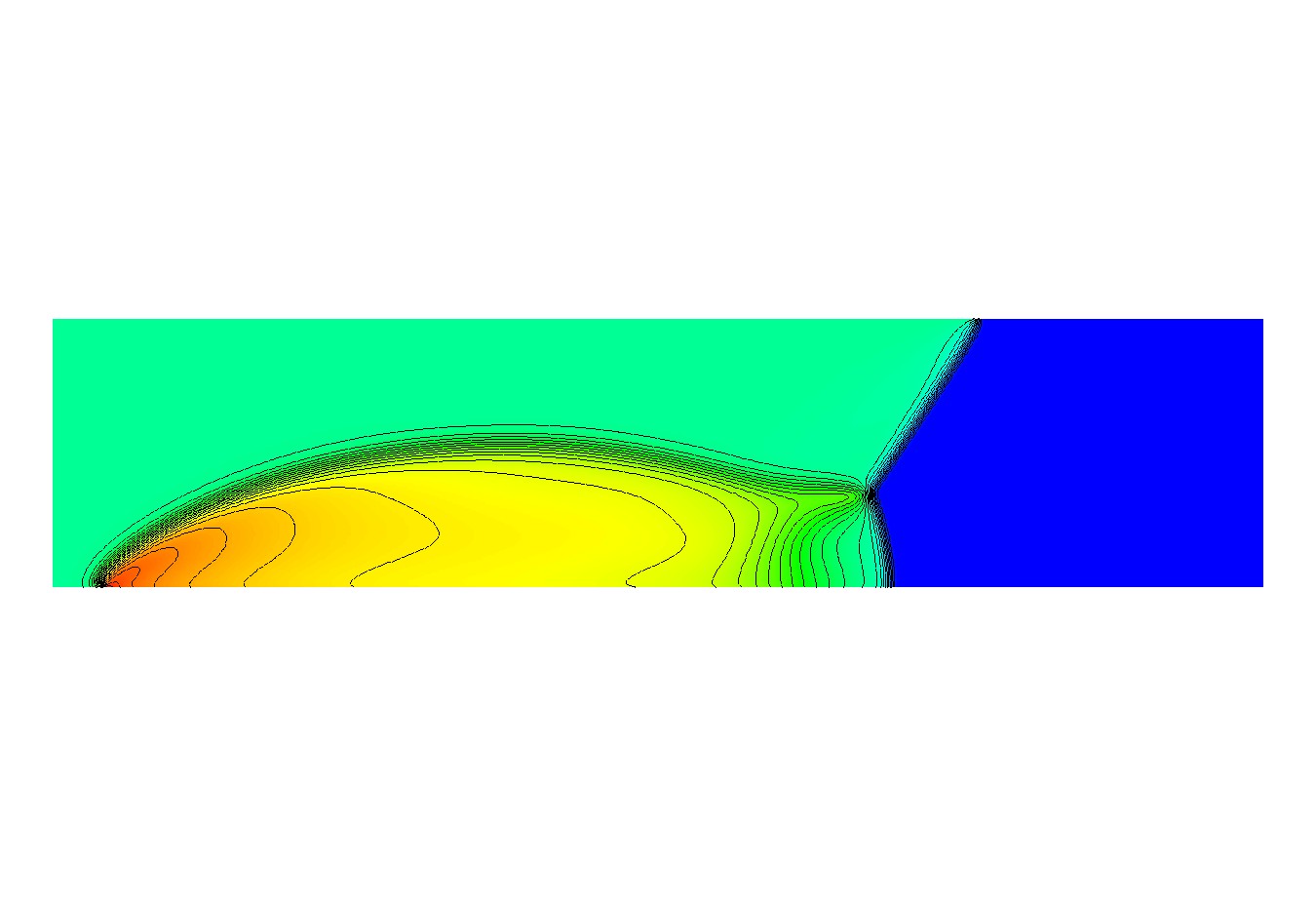}

  (b) MCL-C\&S solution

\includegraphics[width=0.8\textwidth,trim=50 300 250 290,clip]{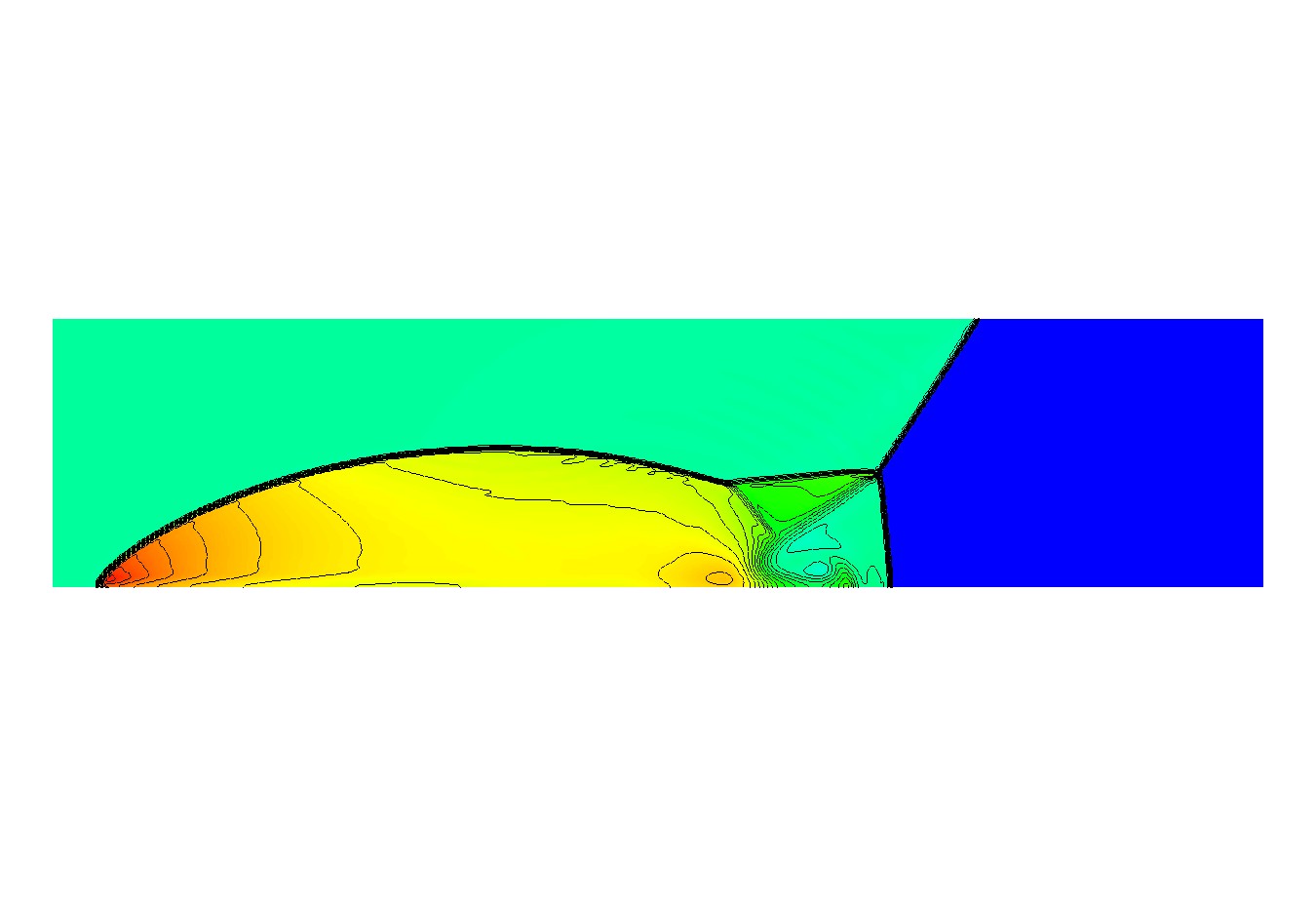}

  \caption{Double Mach reflection: density $\rho_h$ calculated using
    $\mathbb P_1$ elements and $h=1/128$.}
	\label{fig:dmr}      
\end{figure}

Figure \ref{fig:dmr} shows the low-order and MCL-C\&S results for the
double Mach reflection problem \cite{woodward}. A detailed description of the
computational setup can be found, e.g., in \cite{convex}. The presented
numerical solutions are nonoscillatory and satisfy the IDP constraints.
The sequential MCL-C\&S algorithm with IDP correction is clearly
more accurate than the underlying element-based Rusanov method.

\section{Concluding remarks}
\label{sec:5}

The purpose of this review was to highlight the common structure of
element-based convex limiting algorithms and their RD counterparts.
Further examples of such FCT/MCL schemes can be found in \cite[Chapter 4]{book}.
In our experience, recognition of existing relationships
between seemingly different approaches leads to algorithms that
advance the state of the art in all fields. For example, convex
limiting techniques originally designed for discontinuous Galerkin
(DG) methods can be adapted to the continuous finite element setting.
In our recent work \cite{preprint}, we propose an element-based
MCL scheme that constrains intermediate cell averages and nodal states in much
the same way as the Zhang--Shu limiter for the DG version. Moreover,
we stabilize the antidiffusive element contributions using a dissipative
component depending on a weighted essentially nonoscillatory (WENO)
reconstruction. This nonlinear stabilization eliminates
the need to preserve stringent local bounds. Hence, only the
global bounds of IDP constraints need to be enforced using
limiters.

Finally, we remark that the low-order scheme presented
in Section \ref{sec:2} is equivalent to a subface finite volume (FV)
discretization \cite{FV} on a dual mesh $\mathcal T_h^*$ of control volumes
associated with the nodal points $\mathbf x_i$ of a finite element
mesh $\mathcal T_h$. In the FV version, vertices become cell
centers, while cell centers become vertices (cf. \cite{selmin1993,selmin1996}).
This observation makes it possible to incorporate convex limiting of
FCT and MCL type into FV schemes
for general polygonal and polyhedral meshes. The virtual finite element
MCL method proposed by Abgrall et al. \cite{pampa} is also
applicable to such
meshes, and its derivation is based on the RD design philosophy.

\begin{acknowledgement}
  This work is dedicated to Professor Roland Duduchava (University of Georgia) in honor of his 80th anniversary. The author is grateful to him for the opportunity to give a talk on monolithic convex limiting at the Tbilisi Analysis \& PDE online seminar in October 2025.

\end{acknowledgement}


\end{document}